\documentclass[11pt]{article}
\usepackage{graphicx}
\usepackage{amsmath}
\usepackage{amsfonts}
\usepackage{epsfig}
\usepackage{setspace}

\newcommand{\be}{\begin{equation}}
\newcommand{\ee}{\end{equation}}
\newcommand{\beqn}{\begin{eqnarray}}
\newcommand{\eeqn}{\end{eqnarray}}
\newcommand{\beqns}{\begin{eqnarray*}}
\newcommand{\eeqns}{\end{eqnarray*}}

\newcommand{\Var}{\mbox{Var}}
\newcommand{\supp}{\mbox{supp}\ }

\newcommand{\EE}{\ensuremath{{\mathbb E}}}
\newcommand{\II}{\ensuremath{{\mathbb I}}}

\newcommand{\fr}[1]{(\ref{#1})}

\newtheorem{lemma}{Lemma}
\newtheorem{theorem}{Theorem}

\newtheorem{remark}{Remark}

\begin{document}

\title{\Large{\bf  Laplace deconvolution in the presence of indirect  long-memory data }}

\author{
\large{ Rida Benhaddou}  \footnote{E-mail address: Benhaddo@ohio.edu}
  \\ \\
Department of Mathematics, Ohio University, Athens, OH 45701} 
\date{}

\doublespacing
\maketitle
\begin{abstract}
We investigate the problem of estimating a function $f$ based on observations from its noisy convolution when the noise exhibits long-range dependence. We construct an adaptive estimator based on the kernel method, derive minimax lower bound for the  $L^2$-risk when $f$ belongs to Sobolev space and show that such estimator attains optimal rates that deteriorate as the LRD worsens.\\


{\bf AMS (2000) Subject Classification: 62G05, 62G20, 62G08 }
 \end{abstract} 

\section{Introduction.}

Consider the model 

\be
y(t_i)=q(t_i) +\sigma \varepsilon_i ,\ \ \ q(t_i)=\int^{t_i}_0g(t_i-x)f(x)dx, \ \ \ i=1, 2, \cdots, n. \label{conveq}
\ee
where $0\leq t_1 \leq t_2 \leq \cdots \leq t_n \leq T_n$, and the errors $\varepsilon_i$ are Gaussian random variables that are dependent on each other. Let $\varepsilon_n$ be a zero mean vector with components $\varepsilon_i$, $i=1, 2, \cdots, n$, and let $\Sigma_n=Cov(\varepsilon_n)=\EE\left[\varepsilon_n \varepsilon^T_n\right]$ be its covariance matrix. Consider the following assumptions about the errors $\varepsilon_i$ and their covariance matrix $\Sigma_n$.

\noindent
{\bf Assumption A.1.} The vector $\varepsilon_n$ is such that
\be \label{epseta}
\varepsilon_n=A_n\eta_n
\ee
where $\eta_n$ is a vector with independent Gaussian components $\eta_i$, $i=1, 2, \cdots, n$, and $A_n$ is a matrix. Notice that since $\eta_n  \sim \mathcal{N}(0, I_n)$, then under \fr{epseta}, $\Sigma_n=\EE\left[\varepsilon_n \varepsilon^T_n\right]=\EE\left[A_n\eta_n\eta^T_nA^T_n\right]=A_nA^T_n$.

\noindent
{\bf Assumption A.2.}  For the covariance matrix $\Sigma_n$, there exists constants $c_1$ and $c_2$ ($0 < c_1 \leq c_2 < \infty$), independent of $n$, such that
\be
c_1 n^{1-\alpha}\leq \lambda_{min}(\Sigma_n) \leq \lambda_{max}(\Sigma_n)\leq c_2 n^{1-\alpha}, \ \ 0< \alpha \leq 1. \label{LRD}
\ee
where $\alpha$ is the long-memory parameter, and $\lambda_{min}(\Sigma_n)$ and $\lambda_{max}(\Sigma_n)$ are the smallest and the largest eigenvalues of the matrix $\Sigma_n$, respectively.

 Assumption {\bf A.2} is valid when $\varepsilon_n$ is fractional Gaussian or fractional ARIMA, (e.g, see Benhaddou et al.~(2014)). When $\alpha=1$, model \fr{conveq} reduces to the independent and identically distributed noise case.
Model \fr{conveq} with $\alpha=1$ was studied in Abramovich et al.~(2013), where an estimator that relies on the kernel technique is proposed and the choice of the bandwidth is performed by Lepski's Method. Models of this type are referred to as Laplace deconvolution with a noise. This problem is motivated by the analysis of dynamic contrast enhanced imaging data, or modeling time-resolved measurements in fluorescence spectroscopy, (see Abramovich et al.~(2013) for more detail).

Noisy Laplace deconvolution has attracted a lot of attention as of late. One can list a few endeavors such as Dey et al.~(1998), Abramovich et al.~(2013), Vareschi~(2015), and Comte et al.~(2017).  In these attempts, it is assumed that errors  are independent and identically distributed Gaussian random variables. However, empirical evidence has shown that, even at large lags, the correlation structure in the errors can decay at a power-like rate, rather than an exponential rate.   

LRD has been investigated quite considerably in the regression estimation framework, and to some less extent in the standard (Fourier) deconvolution model, (e.g., see Benhaddou et al.~(2014) for more). However, to the best of our knowledge, no LRD work has been published in the context of noisy Laplace deconvolution. 

The objective of the paper is to look into noisy Laplace deconvolution with the relaxation that the noise may exhibit long-range dependence. We do not limit our consideration to one specific form of LRD, rather, we focus our attention to noise structures that satisfy Assumption {\bf A.2}. We establish minimax lower bounds in the $L^2$-risk for estimators of the response function $f$ in model \fr{conveq}  under condition \fr{LRD} when $f$ belongs to a Sobolev ball of radius $A>0$. In addition, we follow the footsteps of Abramovich et al.~(2013) and construct an optimally adaptive estimator that is based on the kernel method, with the optimal choice of the bandwidths performed via Lepski's method. Moreover, we demonstrate that such estimator attains minimax optimal rates. In particular, we show that the convergence rates under LRD depend on a balance between the smoothness parameter of the response function $f$, the parameter of the convolution kernel $g$, and $\alpha$, the long-memory parameter. Finally, it turns out that our convergence rates are similar to those in Abramovich et al.~(2013) when $\alpha=1$, and deteriorate as the level of the LRD gets more and more severe. 

The rest of the paper is organized as follows. Section 2 introduces some notation as well as some assumptions that will be used in the construction of the theoretical results. Section 3 describes the derivation of the lower bounds for the $L^2$-risk of estimators of $f$ observed in model \fr{conveq}. Section 4 goes over the construction of the kernel estimator. Section 5 analyzes the estimation error and shows that when the bandwidths are selected according to Lepski method the estimator attains optimal convergence rates. Finally, Section 6 contains the proofs of the theoretical results.
 \section{Preliminaries.}
 For the rest of the paper, let $\|h\|$ and $\|h\|_{\infty}$ denote the $L^2$-norm and the supremum norm of the function $h$, respectively. Let $\bf{W^{r, p}}$ denote the Sobolev space of the functions defined on $[0, \infty )$ that have $r$ weak derivatives with finite $L^p$-norm, and for $p=2$ denote such space by $\bf{W^r}$. 
 
 Finally, let $r \geq 1$ be such that
 \be \label{Cg}
g^{(j)}(0)= \left\{ \begin{array}{ll} 
 0 , & \mbox{if}\ \  j=0, 1, 2, \cdots, r-2,\\
B_r\neq 0  , & \mbox{if} \ \  j=r-1.
 \end{array} \right.
\ee 
 Next is the list of conditions that will be utilized in the construction of the theoretical results. In particular, the convolution kernel $g$ and the response function $f$ are such that\\
{\bf Assumption A.3.} $g \in \bf{W^{r, 1}}\cap\bf{W^{\nu}}$, with $\nu \geq r$.\\
{\bf Assumption A.4.} Let $\bf{\Omega}$ be the collection of distinct zeros $S_{\omega}$ of the Laplace transform $\tilde{g}$ of $g$. Then, all zeros of $\tilde{g}$ have negative real parts. That is, $$S^*=\max_{S_{\omega\in \Omega}}Re(S_{\omega})<0$$
{\bf Assumption A.5.} $f\in \bf{W^m}$ where $m+r \leq \nu +1$. \\
In addition, $T_n$ satisfies the following\\
{\bf Assumption A.6.} Let $T_n$ be such that $T_n \rightarrow \infty$ but $\frac{T^2_n}{n^{\alpha}}\rightarrow 0$ as $n \rightarrow\infty$, with $0 < \alpha \leq 1$, and there exists $1\leq \mu < \infty$ such that
$$\max_{1 \leq i \leq n}\left|t_i-t_{i-1}\right|\leq \mu \frac{T_n}{n}.$$
\section{Minimax lower bounds.}
In order to establish the performance of estimators for the unknown function $f$ observed in model $(1)$, we derive the minimax lower bounds for the $L^2[0, T_n]$-risk over the Sobolev ball ${\bf{W^m}}(A)$ of radius $A>0$. We define the minimax $L^2$-risk over a set $\Theta$ as 
$$R_n(\Theta)=\inf_{\tilde{f}}\sup_{f\in \Theta}\EE\|\tilde{f}-f\|^2$$
where the infimum is taken over all possible estimators $\tilde{f}$ of $f$. 
The next statement provides the lower bounds of the $L^2[0, T_n]$-risk.
\begin{theorem}\label{th:lowerbds}. Let condition \fr{Cg} and Assumptions {\bf{A.2-A.6}} hold. Then, as $n\rightarrow \infty$,
 \be \label{lowerbds}
 R_n({\bf{W^m}}(A))\geq C\left[ \frac{T_n^2}{n^{\alpha}}\right]^{\frac{2m}{2m +2r +1}}
 \ee
 \end{theorem}
    \section{Estimation Algorithm.}
 In Abramovich et al.~(2013), it was shown that the convolution type Volterra equation of the first kind 
 \be
 q(t)=\int^t_0g(t-x)f(x)dx, \ \ t \geq 0.
 \ee
 admits, under conditions ${\bf{A.3}}$-${\bf{A.5}}$ and \fr{Cg}, the solution 
 \be \label{fsol}
 f(t)=B_r^{-1}\left( q^{(r)}(t)-\sum^{r-1}_{j=0}a_{0, r-j-1}q^{(j)}(t)- \int^t_0q^{(r)}(t-x)\phi_1(x)dx\right)
 \ee
 where  
 \beqn  \label{KLbk:upperbd}
\phi_1(x) &=&\sum^M_{l=1}\sum^{\alpha_l-1}_{j=0} \frac{a_{l, j}X^je^{S_lx}}{j!}\label{fi1}\\
a_{l, j}&=& \frac{1}{(\alpha_l-1-j)!} \frac{d^{\alpha_l-1-j}}{dS^{\alpha_l-1-j}} \left[ (S-S_l)^{\alpha_l}\tilde{\phi}(S)  \right]\label{alj}\\
\tilde{\phi}(S)&=& \frac{S^r\tilde{g}(S)-B_r}{S^r\tilde{g}(S)} \label{fitil}
\eeqn   
and $M$ is the number of distinct zeros of $\tilde{g}(S)$ of orders $\alpha_l$, $l=1, 2, \cdots, M$, and $\alpha_0=r$.

An estimator $\widehat{f}_n$ for $f$ is given by 
\be \label{fsoles}
 \widehat{f}(t)=B_r^{-1}\left( \widehat{q^{(r)}}(t)-\sum^{r-1}_{j=0}a_{0, r-j-1}\widehat{q^{(j)}}(t)- \int^t_0\widehat{q^{(r)}}(t-x)\phi_1(x)dx\right)
 \ee
 where $\widehat{q^{(j)}}(t)$ are some estimators for $q^{(j)}(t)$, $j=0, 1, 2, \cdots, r$. \\
 All we have to do now is use some nonparametric estimation approach to estimate $q(t)=f*g(t)\in \bf{W^{r+m}}$ and its derivatives of orders up to $r$ and plug in \fr{fsoles}. Indeed, we apply the kernel estimation procedure described in Abramovich et al.~(2013). More specifically, choose a kernel function $K_j$ of order $(L, j)$, with $L >r$, that satisfies the following conditions \\
{\bf C.1.} Let $\supp(K_j)=[-1, 1]$, and $\int^1_{-1}K^2_j(t)dt < \infty$. \\
{\bf C.2.} \be \label{mKj}
\int^1_{-1}t^lK_j(t)dt= \left\{ \begin{array}{ll} 
 0 , & \mbox{if}\ \  l=0, 1, 2, \cdots, j-1, j+1, \cdots, L-1,\\
(-1)^jj!  , & \mbox{if} \ \  l=j.
 \end{array} \right.
\ee 
Then, allow the kernel estimator for $q^{(j)}$ with a global bandwidth $\lambda_j$ given by 
\be \label{qjes}
\widehat{q^{(j)}}_{n, \lambda_j}(t)= \frac{1}{\lambda_j^{j+1}}\sum^n_{i=1}K_j\left(\frac{t-t_i}{\lambda_j}\right)(t_i-t_{i-1})y(t_i)
\ee
where the value of $\lambda_j$ is to be determined. 
 \section{Convergence rates and adaptivity.}
 It is necessary to choose the bandwidth levels that minimize the upper bound of the $L^2[0, T_n]$-risk. Indeed, let us investigate the mean integrated squared error. Note that such quantity can be partitioned as follows
 \be \label{eseroll}
\EE \| \widehat{f}_{n}(t)-f(t)\|^2\leq \frac{2+r}{B^2_r}\left[ R_1+ R_2 + R_3\right]
\ee
where 
\beqns
R_1&=& \EE \| \widehat{q^{(r)}}_{n, \lambda_j}(t)-q^{(r)}(t) \|^2\leq 2\int^{T_n}_0\left[\Var\left(\widehat{q^{(r)}}_{n, \lambda_j}(t)\right)+ B^2\left(\widehat{q^{(r)}}_{n, \lambda_j}(t)\right)\right]dt \nonumber\\
R_2&=& \EE \| \left(\widehat{q^{(r)}}_{n, \lambda_j}-q^{(r)}\right)*\phi_1(t) \|^2\leq 2\|\phi_1(t)\|^2\int^{T_n}_0\left[\Var\left(\widehat{q^{(r)}}_{n, \lambda_j}(t)\right)+ B^2\left(\widehat{q^{(r)}}_{n, \lambda_j}(t)\right)\right]dt  \nonumber\\
R_3&=& \sum^{r-1}_{j=0}a^2_{0, r-j-1}\EE \| \widehat{q^{(j)}}_{n, \lambda_j}(t)-q^{(j)}(t) \|^2\leq 2\sum^{r-1}_{j=0}a^2_{0, r-j-1}\int^{T_n}_0\left[\Var\left(\widehat{q^{(j)}}_{n, \lambda_j}(t)\right)+ B^2\left(\widehat{q^{(j)}}_{n, \lambda_j}(t)\right)\right]dt
\eeqns
where $B(\widehat{q^{(j)}}_{n, \lambda_j}(t))$ is the bias of estimator $\widehat{q^{(j)}}_{n, \lambda_j}(t)$. Then, the following statement is true.
\begin{lemma} \label{lem:Var}
Let conditions ${\bf{A.2}}$-${\bf{A.6}}$ hold. Let $\widehat{q^{(j)}}_{n, \lambda_j}(t)$ be defined in \fr{qjes}. Then, for $j=0, 1, \cdots, r$, one has 
\be \label{var-bias}
\sup_{q\in \bf{W^{r+m}}(A')}\EE\|\widehat{q^{(j)}}_{n, \lambda_j}(t)-q^{(j)}(t)\|^2 = O \left (\frac{T^2_n}{n^{\alpha}\lambda^{2j+1}_j}+ \lambda_j^{2(m+r-j)}\right).
\ee
\end{lemma}
Observe that when the bandwidth level $\lambda_j$ increases the bias term increases, while the variance term decreases. It turns out that the estimation error is minimized at the optimal bandwidth level $\lambda_j=\lambda_o$ where both error components are balanced. More specifically, the optimal bandwidth level $\lambda_o$ is given by 
\be  \label{oplam}
 \lambda_o \asymp  \left[ \frac{T_n^{2}}{n^{\alpha}}\right]^{\frac{1}{2(m+r) + 1}}.
\ee
The optimal bandwidth level $\lambda_o$ is expressed in terms of the unknown smoothness parameter $m$ of the function $f$, and therefore it can not be used in the estimation process. One approach to selecting the optimal bandwidth level adaptively is Lepski's method, introduced in Lepski~(1991) and further improved in Lepski et al.~(1997). We borrow some of the ideas of Abramovich et al.~(2013) and adjust them to our setting.

 \noindent\underline{\bf  The Lepski method.} For each $j$, $0 \leq j \leq r$, and the corresponding kernel $K_j$ of order $(L, j)$, $L > r$, consider for some $a>1$ the geometric grid of bandwidths $\Lambda_j$ such that
\be  \label{Lamset}
 \Lambda_j= \left\{  \lambda_j=a^{-l}, l=0, 1, \cdots, J_n: J_n=\frac{1}{2j+1}\log_a\left[ \frac{n^{\alpha}}{\sigma^2T_n^2}\right]\right\}. 
\ee
Define 
\be \label{thres-lev}
\rho_{\lambda, j}^2=\frac{4\sigma^2T_n^2}{n^{\alpha}\lambda_j^{2j+1}}
\ee
Lepski's method suggests to choose the bandwidth level $\lambda_j=\widehat{\lambda}_j$ as 
\be  \label{Lev:Jad}
\widehat{\lambda}_j= \max\left\{  \lambda_j \in  \Lambda_j:  \| \widehat{q^{(j)}}_{n, \lambda_j}(t)-\widehat{q^{(j)}}_{n, \lambda'_j}(t) \|^2\leq \gamma_j^2\rho_{\lambda', j}^2, \ for\ all\  \lambda'_j \in  \Lambda_j, \ \ \lambda'_j < \lambda_j\right\}. 
\ee
In order to see how the selection algorithm performs, observe that 
\be \label{var-op}
\EE\|\widehat{q^{(j)}}_{n, \widehat{\lambda}_j}(t)-q^{(j)}(t)\|^2 = \Delta_1 + \Delta_2.
\ee
where
\beqn 
\Delta_1&=& \EE\left[\|\widehat{q^{(j)}}_{n, \widehat{\lambda}_j}(t)-q^{(j)}(t)\|^2 \II\left(\widehat{\lambda}_j > \lambda_o\right)\right]\label{del1}\\
\Delta_2&=& \EE\left[\|\widehat{q^{(j)}}_{n, \widehat{\lambda}_j}(t)-q^{(j)}(t)\|^2 \II\left(\widehat{\lambda}_j < \lambda_o\right)\right]\label{del2}
\eeqn
If $\widehat{\lambda}_j > \lambda_o$, then by \fr{Lev:Jad}, one has 
\be
\| \widehat{q^{(j)}}_{n, \widehat{\lambda}_j}(t)-\widehat{q^{(j)}}_{n, \lambda_o}(t) \|^2\leq \gamma_j^2\rho_{\lambda_o, j}^2\leq C\left[ \frac{T_n^2}{n^{\alpha}}\right]^{\frac{2(m+r-j)}{2m +2r +1}}
\ee
so that
\be \label{del1up}
\Delta_1\leq C\left[ \frac{T_n^2}{n^{\alpha}}\right]^{\frac{2(m+r-j)}{2m +2r +1}}
\ee
Now, if $\widehat{\lambda}_j < \lambda_o$, then by \fr{Lev:Jad} there exists $\overline{\lambda}_j$,  $\widehat{\lambda}_j < \overline{\lambda}_j< \lambda_o$ such that 
\be
\| \widehat{q^{(j)}}_{n, \overline{\lambda}_j}(t)-\widehat{q^{(j)}}_{n, \lambda_o}(t) \|^2> \gamma_j^2\rho_{\overline{\lambda}_j, j}^2
\ee
It turns out that the probability of such event is very small. In particular, the following lemma provides large deviation results.
\begin{lemma} \label{large-dev}
Let conditions ${\bf{A.1}}$-${\bf{A.6}}$ hold. Let $\widehat{q^{(j)}}_{n, \lambda_j}(t)$be  defined in \fr{qjes}. If the bandwidth is such that 
\be
\left(\frac{\sigma^2T_n^2}{n^{\alpha}}\right)^{\frac{1}{2j+1}}< \lambda< \lambda_o
\ee
and 
\be  \label{Set:Rho}
\gamma_j^2 >\mu \|K_j\|^2
\ee
Then,  as $n \rightarrow \infty$,
\be  \label{Largdev}
\Pr \left(\| \widehat{q^{(j)}}_{n, {\lambda}}(t)-\widehat{q^{(j)}}_{n, \lambda_o}(t) \|^2> \gamma_j^2\rho_{{\lambda}, j}^2\right)= O \left ( \exp\left\{ -\sigma_0(n^{\alpha})^{\frac{1}{2m+2r+1}}(T_n^2)^{\frac{2m +2r -1}{2m +2r +1}}\right\}  \right),
\ee
where $\sigma_0=\frac{\|K_j\|^2}{2c_2}$ and $c_2$ appears in condition \fr{LRD} and $\mu$ appears in Assumption ${\bf{A.6}}$. 
\end{lemma}
Consequently, the following theorem provides the upper bounds of $L^2[0, T_n]$-risk of the estimation based on a bandwidth selection according to \fr{Lev:Jad}. 
\begin{theorem} \label{th:upperbds}
Let conditions ${\bf{A.1}}$-${\bf{A.6}}$ and \fr{Cg} hold. Let $\tilde{f}_{n}(t)$ be defined in \fr{fsoles} with $\widehat{g^{(j)}}_{n, \lambda_j}(t)$ given in \fr{qjes}. Choose the bandwidth $\lambda_j=\widehat{\lambda}_j$ according to \fr{Lev:Jad} with $\gamma_j^2$ satisfying \fr{Set:Rho}. Then, for all $1\leq m\leq \min(L, \nu+1)-r$, as $n \rightarrow \infty$, one has
 \be \label{upperbds}
 R_n({\bf{W^m}}(A))\leq C\left[ \frac{T_n^2}{n^{\alpha}}\right]^{\frac{2m}{2m +2r +1}}
 \ee
\end{theorem}

 \begin{remark}	
{\rm{
{\bf{(i)}}\,
Note that under additional but minor conditions on $f$ and $T_n$, Theorems \ref{th:lowerbds} and \ref{th:upperbds} can be extended to $t \in [0, \infty )$.  \\
 \noindent {\bf{(ii)}}\,
Theorems \ref{th:lowerbds} and \ref{th:upperbds} imply that, for the $L^2$-risk, the estimator \fr{fsoles} with $\widehat{g^{(j)}}_{n, \lambda_j}(t)$ given by \fr{qjes} and global bandwidths $\lambda_j=\widehat{\lambda}_j$ according to \fr{Lev:Jad} is adaptive and asymptotically optimal over all  Sobolev spaces  $ W^{m}(A)$. \\  
  \noindent {\bf{(iii)}}\, 
The convergence rates are expressed in terms of the long-memory parameter $\alpha$, in addition to the parameters $m$ and $r$ associated with the smoothness of the functions $f$ and $g$, respectively. In particular, the rates deteriorate as the long-range dependence gets more severe. This behavior is consistent with that in Benhaddou~(2016), Kulik et al.~(2015), Benhaddou et al.~(2014) or Wishart~(2013), in their standard (Fourier) deconvolution with LRD setup.\\   \noindent {\bf{(iv)}}\,
 For  $\alpha=1$ our rates match exactly those in  Abramovich et al.~(2013) in their case of Laplace deconvolution with i.i.d. noise.
   }}
\end{remark}

 \section{Proofs. }

  {\bf Proof of Theorem \ref{th:lowerbds}}.
In order to prove the theorem, we consider the test functions used in the construction of lower bounds obtained by Abramovich et al.~(2013). Lemma $A.1$ of Bunea et al.~(2007) is then applied to find such lower bounds using Assumptions  {\bf{A.2-A.6}}, along with condition \fr{Cg}. Assume that the points $t_i$ are equally spaced. Define the integers $M_n \geq 8$ and $N=\left[\frac{n}{M_n}\right]$, the largest integer which does not exceed $\frac{n}{M_n}$. Let $\lambda_n=N\frac{ T_n}{n}$ and define the points $z_j=j\lambda_n$, $j=0, 1, \cdots, M_n$. Keep in mind that based on the definition of $N$, $\frac{T_n}{2M_n} \leq \lambda_n\leq \frac{T_n}{M_n}$. Let $k(.)$ be infinitely differentiable function, with $\supp(k)=[0, 1]$, such that
\be
\int^1_0x^jk(x)dx=0,\ \ j=0, 1\cdots, r-1, \ \ and\ \ \int^1_0x^rk(x)dx\neq 0.
\ee
and introduce the functions
\be
\psi_j(x)=L \frac{\lambda_n^m}{\sqrt{T_n}}k\left(\frac{x-z_{j-1}}{\lambda_n}\right), \ \ j=1, 2, \cdots, M_n.
\ee
Let $\omega$ be a vector with components $\omega_j=\{0, 1\}$, $j=1, 2, \cdots, M_n$. Denote the set of all possible values of $\omega$ by $\Theta$ and let the functions $f_{\omega}$ be of the form
\be  \label{f_j}
f_{\omega}(t)=\sum^{M_n}_{j=1}\omega_j\psi_j(t), \ \ \omega\in \Theta.
\ee
It is easy to verify that  $f_{\omega}(t) \in W^{m}(A)$ with $A=L\|k\|$. If $f_{\tilde{\omega}}$ is the form \fr{f_j} but with $\tilde{\omega} \in \Theta$ instead of $\omega$, then the $L^2$-norm of the difference is
\be
\| f_{\omega}(t)-f_{\tilde{\omega}}(t)\|^2=L^2\frac{\lambda_n^{2m+1}}{T_n}H(\tilde{\omega}, \omega)
\ee
where $H(\tilde{\omega}, \omega)$ is the Hamming distance between the binary sequences $\omega$ and $\tilde{\omega}$. 
Remark that vector $\omega$ has ${M_n}$ components, and therefore, $Card(\Theta)=2^{M_n}$. In order to find a lower bound for $H(\tilde{\omega}, \omega)$ we apply the Varshamov-Gilbert Lemma which argues that one can choose a subset $\Theta_1$ of $\Theta$, of cardinality of at least $2^{M_n/8}$ such that for any $\omega,\ \tilde{\omega}\in\Theta_1$, $H(\tilde{\omega}, \omega)\geq \frac{M_n}{8}$. Hence,
\be \label{normdif}
\| f_{\omega}(t)-f_{\tilde{\omega}}\|^2 \geq \frac{L^2\|k(x)\|^2}{8}\lambda_n^{2m}=4\delta^2
\ee
Now define $q^{(n)}_{\omega}$, the vector with components
\be
q_{\omega}(t_i)=\int^{t_i}_0g(t_i-x)f_{\omega}(x)dx, \ \ i=1, 2, \cdots,n.
\ee
In addition, define the quantities 
\be
k_1(x)=\int^x_0k(t)dt,\ \  \ k_j(x)=\int^x_0k(t)dt,\ \ j=2,\cdots, r.
\ee
Then, using the fact that $|\omega_j-\tilde{\omega}_j|\leq 1$, Lemma 2, results $(5.8)$ and $(5.9)$ in Abramovich et al.~(2013) and Assumption ${\bf A.2}$, the Kullback divergence can be written as 
\beqn
K(P_{f_{\omega}}, P_{f_{\tilde{\omega}}})&=&\frac{1}{2\sigma^2}\left(q^{(n)}_{\omega}- q^{(n)}_{\tilde{\omega}}\right)^T(\Sigma_n)^{-1} \left(q^{(n)}_{\omega}- q^{(n)}_{\tilde{\omega}}\right)\nonumber\\
&\leq& \frac{1}{2\sigma^2}\lambda_{\max}\left(\Sigma_n\right)^{-1}\sum^n_{i=1}\left[q_{\omega}(t_i)- q_{\tilde{\omega}}(t_i)\right]^2\nonumber\\
&\leq& \frac{L^2\lambda_n^{2m+2r}}{\sigma^2c_1T_n}n^{\alpha-1}(E_1+ E_2)
\eeqn
where 
\be
E_1= \sum^n_{i=1}\left[ \sum^{M_n}_{j=1}B_rk_r\left(\frac{t_i-z_{j-1}}{\lambda_n}\right)\II\left(z_{j-1} \leq y_i\leq z_j\right)\right]^2\leq nB_r^2\|k_r(x)\|^2_{\infty}
\ee
and
\be
E_2=\sum^n_{i=1}\left[\sum^{M_n}_{j=1}\int^{\min\{z_j, t_i\}}_{\min\{z_{j-1}, t_i\}}g^{(r)}(t_i-x)k_r\left(\frac{x-z_{j-1}}{\lambda_n}\right)dx\right]^2\leq n\|g^{(r)}\|^2\|k_r(x)\|^2_{\infty}
\ee
Consequently, the application of Lemma $A.1$  requires 
\be \label{klblam}
L^2\lambda_n^{2(m+r)}\frac{n^{\alpha}}{\sigma^2c_1T_n}\leq CM_n
\ee
and, since $\frac{T_n}{2M_n}\leq \lambda_n\leq \frac{T_n}{M_n}$, one obtains 
\be \label{mn}
\left(T_n\right)^{2m +2r-1}n^{\alpha}\leq C\left(M_n\right)^{2m +2r+1}
\ee
It is easy to check that  the choice 
\be \label{lamn}
\lambda_n=\left(\frac{T_n^2}{n^{\alpha}}\right)^{\frac{1}{2m+2r+1}}
\ee
satisfies \fr{klblam}. Hence, by Lemma A.1 we conclude the lower bound by plugging \fr{lamn} in \fr{normdif}. $\Box$
{\bf Proof of Lemma \ref{lem:Var}.} Observe that the quantities $\aleph_{j, \lambda_j}(t)=\widehat{q^{(j)}}_{n, {\lambda}_j}(t) - \EE\left[\widehat{q^{(j)}}_{n, {\lambda}_j}(t)\right]$ are zero mean Gaussian random variables with variance 
\be
\Var\left[\aleph_{j, \lambda_j}(t)\right]=\frac{\sigma^2}{\lambda_j^{2j+2}}K_{n, j}^T(t)\Sigma_nK_{n, j}(t)
\ee
where $K_{n, j}(t)$ are vectors with elements $K_j\left(\frac{t-t_i}{\lambda_j}\right)(t_i-t_{i-1})$, $i=1, 2, \cdots, n$, and $\Sigma_n$ is the covariance matrix of the vector $\varepsilon_n$. Therefore, by Assumption ${\bf{A.2}}$, the variance becomes
\beqn
\Var\left[\aleph_{j, \lambda_j}(t)\right]&\leq& \frac{\sigma^2}{\lambda_j^{2j+2}}\lambda_{\max}\left[\Sigma_n\right]\|K_{n, j}(t)\|^2\nonumber\\
&\leq& \frac{2\sigma^2c_2}{\lambda_j^{2j+1}}\frac{T_n}{n^{\alpha}}\|K_j\|^2\label{varalep}
\eeqn
Consequently, integrating both sides of \fr{varalep} over the interval $[0, T_n]$ yields the first term in \fr{var-bias}.

For the bias term, the derivation will be the same as in Abramovich et al.~(2013), so we skip it. This completes the proof. $\Box$\\
{\bf Proof of Lemma \ref{large-dev}.} Following Abramovich et al.~(2013), denote 
\be
d_j=\frac{\gamma_j-\mu \|K_j\|}{2\|K_j\|}
\ee
 and set 
\be
\lambda_o=\left(\frac{d_j^2c_2\sigma^2B_0T_n^2}{2(A')^2n^{\alpha}}\right)^{\frac{1}{2(m+r)+1}}
\ee
where $B_0$ appears in the calculation of the bias term in Abramovich et al.~(2013). Since, $\lambda<\lambda_o$, then, by \fr{var-bias} and \fr{oplam} the bias term is such that
\be
\|B(\widehat{q^{(j)}}_{n, \lambda}(t))\|^2\leq d_j^2c_2\sigma^2\|K_j\|^2\frac{T_n^2}{\lambda^{2j+1}n^{\alpha}}
\ee
Hence,
\beqn
\Pr\left(\|\tilde{q^{(j)}_{\lambda_o}}-\tilde{q^{(j)}_{\lambda}}\|^2> \gamma_j^2\rho^2_{\lambda}\right)&\leq& \Pr\left(\|\aleph_{j, \lambda}(t)\|^2> c_2\sigma^2\|K_j\|^2(\mu + d_j)^2\frac{T_n^2}{n^{\alpha}\lambda^{2j+1}}\right)\nonumber\\
&+& \Pr\left(\|\aleph_{j, \lambda_o}(t)\|^2> c_2\sigma^2\|K_j\|^2(\mu + d_j)^2\frac{T_n^2}{n^{\alpha}\lambda^{2j+1}}\right)\nonumber\\
&\leq& \Pr\left( \varepsilon_n^TQ_{\lambda}\varepsilon_n> c_2\sigma^2\|K_j\|^2(\mu +d_j)^2\right)\nonumber\\
&+& \Pr\left( \varepsilon_n^TQ_{\lambda_o}\varepsilon_n> c_2\sigma^2\|K_j\|^2(\mu +d_j)^2\right)
\eeqn
where $\aleph_{j, \lambda}(t)=\widehat{q^{(j)}}_{n, {\lambda}}(t) - \EE\left[\widehat{q^{(j)}}_{n, {\lambda}}(t)\right]$, and $Q_{\lambda}$ is a symmetric nonnegative definite matrix with elements
\be
Q_{il, \lambda}=\frac{n^2}{T_n^2}(t_i-t_{i-1})(t_l-t_{l-1})\int^1_{-1}K_j(z)K_j(z + \frac{t_i-t_l}{\lambda})dt
\ee
Now we use large deviation result that was developed in Comte.~(2001) and further improved in Gendre~(2014)  which states that for any $x>0$, if $\xi_n$ is a zero mean Gaussian vector with independent elements, and $Q$ is nonnegative definite matrix, then 
\be \label{QLarge-D}
\Pr\left(\xi_n^TQ\xi_n> \sigma^2\left[\sqrt{(tr(Q))}+\sqrt{x\rho^2_{\max}(Q)}\right]\right)\leq e^{-x}
\ee
Therefore, by Assumption {\bf A.1}, the vector $\varepsilon_n$ allows the representation $\varepsilon_n=A_n\eta_n$, where $\eta_n$ is a zero mean Gaussian vector with independent elements, and $A_n$ is a matrix. Hence, $\varepsilon_n^TQ_{\lambda}\varepsilon_n=\eta_n^T[A_n^TQ_{\lambda}A_n]\eta_n$, where $ A_n^TQ_{\lambda}A_n$ is nonnegative definite matrix. In addition, the covariance of $\varepsilon_n$ is $\Sigma_n=A_nA_n^T$. To apply result \fr{QLarge-D}, all we have to do is find $Tr\left(A_n^TQ_{\lambda}A_n\right)$ and $\rho_{\max}^2\left(A_n^TQ_{\lambda}A_n\right)$ and choose an appropriate $x$. Indeed, 
\be
Tr\left(A_n^TQ_{\lambda}A_n\right)=Tr\left(Q_{\lambda}A_nA_n^T\right)\leq \lambda_{\max}\left(\Sigma_n\right)Tr(Q_{\lambda})\leq c_2n^{1-\alpha}n\mu^2\|K_j\|^2
\ee
and
\be
\rho^2_{\max}\left(A_n^TQ_{\lambda}A_n\right)\leq \rho^2_{\max}\left(\Sigma_n\right)\rho^2_{\max}\left(Q_{\lambda}\right)\leq c_2n^{1-\alpha}2\nu \|K_j\|^2\frac{n \lambda}{T_n}
\ee
Finally, applying result \fr{QLarge-D} with the choice $x=d_j^2\frac{T_n}{2c_2\mu \lambda}$, for $\lambda < \lambda_o$, completes the proof.  $\Box$
 {\bf Proof of Theorem \ref{th:upperbds}}. In order to find the upper bound for $\Delta_2$ in \fr{del2}, note that 
\be
\|B(\widehat{q^{(j)}}_{n, \lambda}(t))\|^4=o\left( \lambda^{4(m+r-j)}\right)
\ee
and for any, ${\lambda}< \lambda_o$, then
\be
\EE\|\widehat{q^{(j)}}_{n, \lambda}(t)-q^{(j)}(t)\|^4=O \left(\left(\frac{T^2_n}{n^{\alpha}\lambda^{2j+1}}\right)^2+ \lambda^{4(m+r-j)}\right)=O(1)
\ee
Consequently, by Lemma \ref{large-dev} and Cauchy-Schwarz inequality, $\Delta_2$ in \fr{del2} is such that 
\beqn
\Delta_2 &\leq& \sum^{\lambda_o-1}_{\lambda= \left(\frac{T^2_n}{n^{\alpha}}\right)^{\frac{1}{2j+1}}}\sqrt{\EE\|\widehat{q^{(j)}}_{n, \lambda}(t)-q^{(j)}(t)\|^4}\sqrt{\Pr\left(\| \widehat{q^{(j)}}_{n, {\lambda}}(t)-\widehat{q^{(j)}}_{n, \lambda_o}(t) \|^2> \gamma_j^2\rho_{{\lambda}, j}^2\right)}\nonumber\\
&=&o\left(\left[ \frac{T_n^2}{n^{\alpha}}\right]^{\frac{2(m+r-j)}{2m +2r +1}}\right)\label{Del}
\eeqn
Hence, combining \fr{del1up} and \fr{Del} completes the proof. $\Box$


\begin{thebibliography}{99}

\bibitem{abram} Abramovich, F., Pensky, M., Rozenholc, Y. (2013). Laplace deconvolution with noisy observations.
{\it Electron. J. Stat.} {\bf 7}, 1094-1128.
 \bibitem{ben1} Benhaddou, R. (2016). Deconvolution model with fractional Gaussian noise: a minimax study. {\it Statistics and Probability Letters}. {\bf 117}, 201-208.
 
  \bibitem{ben1} Benhaddou, R., Kulik, R., Pensky, M., Sapatinas, T. (2014). Multichannel deconvolution with long-range dependence: a minimax study. {\it Invited paper at Journ. Statist. Plan. Inf.} {\bf 148}, 1-19.

\bibitem{bun} Bunea, F., Tsybakov, A., Wegkamp, M.H. (2007). Aggregation for Gaussian regression. {\it Ann. Statist.} {\bf 35}, 1674-1697.

\bibitem{com1}
Comte, F. (2001).
Adaptive estimation of the spectrum of a stationary Gaussian sequence.
{\it Bernoulli.} {\bf 7},   267--298.

\bibitem{com2}
Comte, F.,    Cuenod, C.-A.,    Pensky, M.,   Rozenholc, Y. (2017).
Laplace deconvolution on the basis of   time domain data  
and its application to Dynamic Contrast Enhanced imaging.
{\it Journ. Royal Stat. Soc., Ser.B.} {\bf 79},   69--94.

\bibitem{dey}
Dey, A.K., Martin, C.F.,  Ruymgaart, F.H.  (1998). 
Input recovery from noisy output data, using regularized inversion of
Laplace transform. {\it IEEE Trans. Inform. Theory.} {\bf 44}, 1125--1130.

\bibitem{gen} Gendre, X. (2014). Model selection and estimation of a component in additive  regression. {\it ESAIM: Probability and Statistics.} {\bf 18}, 77-116.


\bibitem{kul4} Kulik, R., Sapatinas, T. Wishart, J. R. (2015). Multichannel deconvolution with long-range dependence: Upper bounds on the $L_p$-risk. {\it Applied and Computational Harmonic Analysis.} {\bf 38}, 357-384.

\bibitem{lep1} Lepski, O. V. (1991). Asymptotic minimax adaptive estimation I: Upper bounds. Optimally adaptive estimates. {\it Theory Probab. Appl.} {\bf 36}, 654-659.

\bibitem{lep2} Lepski, O. V., Mammen, E., Spokoiny, V. G. (1991). Optimal spatial adaptation to inhomogeneous smoothness: an approach based on kernel estimates with variable bandwidth selectors. {\it Ann. Statist.} {\bf 25}, 629-947.

 \bibitem{tsybakov}
Tsybakov, A.B. (2008).
{\it Introduction to Nonparametric Estimation}, Springer, New York.

\bibitem{vareschi}
Vareschi, T. (2015).
Noisy Laplace deconvolution with error in the operator.  
{\it Journ. Statist. Plan. Inf.} {\bf 157-158},  16-35.

\bibitem{wan1} Wang, Y. (1997). Minimax estimation via wavelets for indirect long-memory data. {\it Journ. Statist. Plan. Inf.} {\bf 1}, 45-55.

\bibitem{wish} Wishart, J. M. (2013). Wavelet deconvolution in a periodic setting with long-range Dependent  Errors. {\it Journ. Statist. Plan. Inf.} {\bf 5}, 867-881.


\end{thebibliography}
\end{document}